\documentclass[letterpaper,conference]{ieeeconf}
\IEEEoverridecommandlockouts                              
\usepackage[utf8]{inputenc}
\usepackage{textcomp}

\usepackage{graphicx}
\usepackage{amsmath}
\usepackage[version=4]{mhchem}
\usepackage{siunitx}
\usepackage{graphicx}
\usepackage{longtable,tabularx}
\usepackage[style=ieee]{biblatex}
\usepackage{tikz}
\usetikzlibrary{positioning, arrows.meta}
\usepackage{amsmath} 
\usepackage{amsfonts}
\usepackage{amsthm}

\newtheorem{theorem}{Theorem}

\theoremstyle{definition}

\newtheorem{remark}{Remark}

\usepackage{float} 
\usepackage[skip=0em,font=footnotesize]{caption}

\usepackage[margin = 1in]{geometry}

\usepackage{tikz}
\usetikzlibrary{shapes,arrows,calc,positioning}
\tikzstyle{bigblock} = [draw, fill=blue!20, rectangle, 
    minimum height=6em, minimum width=8em]
\tikzstyle{medblock} = [draw, fill=blue!20, rectangle, 
    minimum height=4em, minimum width=4em]    
\tikzstyle{mux} = [draw, fill=black!20, rectangle, 
    minimum height=5em, minimum width=0.1em]    
\tikzstyle{smallblock} = [draw, fill=blue!20, rectangle, 
    minimum height=2em, minimum width=3em]
    
\tikzstyle{data_block} = [draw, fill=green!20, rectangle, 
    minimum height=2em, minimum width=3em]
\tikzstyle{ops_block} = [draw, fill=blue!20, rectangle, 
    minimum height=2em, minimum width=3em]    
\tikzstyle{est_block} = [draw, fill=red!20, rectangle, 
    minimum height=2em, minimum width=3em]    
    
\tikzstyle{sum} = [draw, fill=blue!20, circle, node distance=1cm,minimum height=0.5cm]
\tikzstyle{signal} = [coordinate]
\tikzstyle{pinstyle} = [pin edge={to-,thin,black}]
\tikzstyle{block} = [draw, fill=blue!20, rectangle, 
    minimum height=3em, minimum width=9em]
\tikzstyle{blockS} = [draw, fill=blue!20, rectangle, 
    minimum height=3em, minimum width=4em]    
\tikzstyle{input} = [coordinate]
\tikzstyle{output} = [coordinate]
\usetikzlibrary{matrix}

\usetikzlibrary{positioning}
\usetikzlibrary{math}

\usepackage{hyperref}
\usepackage{xcolor}
\hypersetup{
    colorlinks,
    linkcolor={blue!100!black},
    citecolor={blue!50!black},
    urlcolor={blue!80!black}
}

\newcommand{\bc}{\begin{center}}
\newcommand{\ec}{\end{center}}
\newcommand{\benum}{\begin{enumerate}}
\newcommand{\eenum}{\end{enumerate}}
\newcommand{\nn}{\nonumber}
\newcommand{\matl}{\left[ \begin{array}}
\newcommand{\matr}{\end{array} \right]}

\renewcommand{\matl}{\begin{bmatrix}}
\renewcommand{\matr}{\end{bmatrix}}

\newcommand{\matls}{\left[ \begin{smallmatrix}}
\newcommand{\matrs}{\end{smallmatrix} \right]}
\newcommand{\isdef}{\stackrel{\triangle}{=}}

\newcommand{\vect}[1]{\overset{\rightharpoonup}{#1}}

\newcommand{\rmA}{{\rm A}}

\newcommand{\rmT}{{\rm T}}

\newcommand{\rmc}{{\rm c}}
\newcommand{\rmd}{{\rm d}}

\newcommand{\BBR}{{\mathbb R}}

\newcommand{\SG}{{\mathcal G}}

\title{Backstepping Control of a Bicopter with Unknown Mass and Inertia}

\title{Globally Asymptotically Stable Adaptive Control of a \\ 
Bicopter with Unknown Mass and Inertia}

\title{Estimation inversion-free Adaptive Control of a \\ 
Bicopter with Unknown Mass and Inertia}

\title{Singularity-free Backstepping-based Adaptive Control of a \\ 
Bicopter with Unknown Mass and Inertia}

\author{
Jhon Manuel Portella Delgado
and
Ankit Goel
\thanks{Jhon Manuel Portella Delgado is a graduate student in the Department of Mechanical Engineering, University of Maryland, Baltimore County, 1000 Hilltop Circle, Baltimore, MD 21250. {\tt\small jportella@umbc.edu}}%
\thanks{Ankit Goel is an Assistant Professor in the Department of Mechanical Engineering, University of Maryland, Baltimore County,1000 Hilltop Circle, Baltimore, MD 21250. {\tt\small ankgoel@umbc.edu }}%
}

\begin{document}

\maketitle


\begin{abstract}

The paper develops a singularity-free backstepping-based adaptive control for stabilizing and tracking the trajectory of a bicopter system.
In the bicopter system, the inertial parameters parameterize the input map.
Since the classical adaptive backstepping technique requires the inversion of the input map, which contains the estimate of parameter estimates, the stability of the closed-loop system cannot be guaranteed due to the inversion of parameter estimates.
This paper proposes a novel technique to circumvent the inversion of parameter estimates in the control law. 
The resulting controller requires only the sign of the unknown parameters. 
The proposed controller is validated in simulation for a smooth and nonsmooth trajectory-tracking problem. 

\end{abstract}
\textit{\bf keywords:} adaptive backstepping control, dynamic extension, bicopter, multicopter.

\section{INTRODUCTION}

Multicopters are extensively used in applications like precision agriculture \cite{mukherjee2019}, environmental surveys \cite{lucieer2014,klemas2015}, construction management \cite{li2019}, and load transportation \cite{villa2020}. Their growing popularity drives interest in enhancing their capabilities, but reliable control remains challenging due to nonlinear, time-varying dynamics and uncertain environments.

Common control techniques for multicopters \cite{nascimento2019,marshall2021,castillo2004} depend on accurate model parameters, making them sensitive to uncertainties \cite{emran2018,amin2016}. Adaptive methods like model reference adaptive control \cite{whitehead2010,dydek2012}, L1 adaptive control \cite{zuo2014}, adaptive sliding mode control \cite{espinoza2021trajectory,wu20221, mofid2018}, and retrospective cost adaptive control \cite{goel_adaptive_pid_2021,spencer2022} address some of these issues, but often lack stability guarantees or require pre-existing controllers.
Conventional control architectures decompose multicopter dynamics into separate translational and rotational loops \cite{px4_architecture}, simplifying design but not guaranteeing closed-loop stability.
These approaches rely on the time separation principle, where successive loops are designed to be faster than the preceding loops.

This paper develops an adaptive controller for fully coupled nonlinear dynamics using a bicopter system, a simplified multicopter retaining key nonlinear dynamics. 
In particular, we develop an adaptive backstepping controller \cite{krstic1995nonlinear}, overcoming the non-invertibility of the input map via dynamic extension \cite{descusse1985decoupling}. 
A bicopter is a multicopter constrained to a vertical plane and thus is modeled by a 6th-order nonlinear instead of a 12th-order nonlinear system. 
Despite the lower dimension of the state space, the 6th-order bicopter retains the complexities of the nonlinear dynamics of an unconstrained multicopter. 

As shown in our previous work reported in \cite{portella2024adaptive}, the classical backstepping technique can not be used due to the non-invertibility of the input map.
In a bicopter system, the input map is a tall matrix, and thus, the non-invertibility is structural rather than due to rank deficiency at a few points in the state space. 
As shown in \cite{portella2024adaptive}, a dynamic extension of the equations of motion overcomes the structural non-invertibility of the input map. However, the new input map is still parameterized by the unknown parameters.

Unlike existing methods \cite{zhang2017dynamic,mazenc2018backstepping,reger2019dynamic,triska2021dynamic}, which apply to the systems with uncertainties in the dynamics map, the proposed technique is developed for systems with uncertainties in the input map.
In the proposed method, the requirement of inversion of the input map is relaxed to the knowledge of the sign of the unknown parameters. 
This requirement is often benign since the sign of unknown parameters is known due to their physical interpretation. 
For example, inertial properties such as mass and inertia, which may be unknown, are positive. 
The key contribution of this paper is thus the extension of the backstepping framework to design an adaptive controller for a system with uncertainty in the input map, its application to design a control law with stability guarantee for a bicopter system, and validation in simulation.



The paper is organized as follows. 
Section \ref{sec:prob_formulation} presents the equations of motion of the bicopter system.
Section \ref{sec:adaptive_backstepping} describes the adaptive backstepping procedure to design the adaptive controller for the bicopter.
Section \ref{sec:simulations} presents the results of numerical simulations to validate the adaptive controller.
Finally, the paper concludes with a discussion of results and future research directions in section 
\ref{sec:conclusions}.

\section{Problem formulation}
\label{sec:prob_formulation}

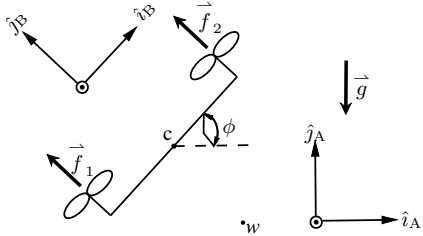
\begin{figure}[!ht]
    \centering
    \resizebox{0.75\columnwidth}{!}{
    
\tikzset{every picture/.style={line width=0.75pt}} 

\begin{tikzpicture}[x=0.75pt,y=0.75pt,yscale=-0.7,xscale=0.7]

\draw   (382.73,198.06) .. controls (384.44,196.15) and (387.38,195.98) .. (389.3,197.69) .. controls (391.22,199.41) and (391.38,202.35) .. (389.67,204.26) .. controls (387.96,206.18) and (385.02,206.35) .. (383.1,204.63) .. controls (381.18,202.92) and (381.02,199.98) .. (382.73,198.06) -- cycle ;
\draw    (385.37,137.84) -- (385.88,195.97) ;
\draw [shift={(385.35,135.84)}, rotate = 89.5] [fill={rgb, 255:red, 0; green, 0; blue, 0 }  ][line width=0.08]  [draw opacity=0] (12,-3) -- (0,0) -- (12,3) -- cycle    ;
\draw    (449.36,199.25) -- (391.23,200.02) ;
\draw [shift={(451.36,199.23)}, rotate = 179.24] [fill={rgb, 255:red, 0; green, 0; blue, 0 }  ][line width=0.08]  [draw opacity=0] (12,-3) -- (0,0) -- (12,3) -- cycle    ;
\draw  [fill={rgb, 255:red, 0; green, 0; blue, 0 }  ,fill opacity=1 ] (384.93,201.16) .. controls (384.93,200.46) and (385.5,199.9) .. (386.2,199.9) .. controls (386.9,199.9) and (387.47,200.46) .. (387.47,201.16) .. controls (387.47,201.86) and (386.9,202.43) .. (386.2,202.43) .. controls (385.5,202.43) and (384.93,201.86) .. (384.93,201.16) -- cycle ;
\draw    (322.35,83.75) -- (220.35,195.75) ;
\draw [shift={(271.35,139.75)}, rotate = 132.32] [color={rgb, 255:red, 0; green, 0; blue, 0 }  ][fill={rgb, 255:red, 0; green, 0; blue, 0 }  ][line width=0.75]      (0, 0) circle [x radius= 1.34, y radius= 1.34]   ;
\draw    (200.75,177.75) -- (220.35,195.75) ;
\draw   (183.73,197.05) .. controls (181.5,194.8) and (183.84,188.88) .. (188.95,183.81) .. controls (194.05,178.75) and (200,176.46) .. (202.22,178.71) .. controls (204.45,180.95) and (202.11,186.88) .. (197.01,191.94) .. controls (191.9,197.01) and (185.95,199.29) .. (183.73,197.05) -- cycle ;
\draw   (202.22,178.71) .. controls (200,176.46) and (202.33,170.54) .. (207.44,165.47) .. controls (212.55,160.41) and (218.49,158.12) .. (220.72,160.37) .. controls (222.94,162.61) and (220.61,168.54) .. (215.5,173.6) .. controls (210.39,178.66) and (204.45,180.95) .. (202.22,178.71) -- cycle ;
\draw    (302.75,65.75) -- (322.35,83.75) ;
\draw   (285.73,85.05) .. controls (283.5,82.8) and (285.84,76.88) .. (290.95,71.81) .. controls (296.05,66.75) and (302,64.46) .. (304.22,66.71) .. controls (306.45,68.95) and (304.11,74.88) .. (299.01,79.94) .. controls (293.9,85.01) and (287.95,87.29) .. (285.73,85.05) -- cycle ;
\draw   (304.22,66.71) .. controls (302,64.46) and (304.33,58.54) .. (309.44,53.47) .. controls (314.55,48.41) and (320.49,46.12) .. (322.72,48.37) .. controls (324.94,50.61) and (322.61,56.54) .. (317.5,61.6) .. controls (312.39,66.66) and (306.45,68.95) .. (304.22,66.71) -- cycle ;
\draw  [dash pattern={on 4.5pt off 4.5pt}]  (331.49,138.96) -- (271.35,139.75) ;
\draw  [fill={rgb, 255:red, 0; green, 0; blue, 0 }  ,fill opacity=1 ] (326.07,201.55) .. controls (326.07,200.98) and (326.53,200.51) .. (327.1,200.51) .. controls (327.68,200.51) and (328.14,200.98) .. (328.14,201.55) .. controls (328.14,202.12) and (327.68,202.59) .. (327.1,202.59) .. controls (326.53,202.59) and (326.07,202.12) .. (326.07,201.55) -- cycle ;
\draw [line width=1.5]    (410.42,70.59) -- (410.12,109.92) ;
\draw [shift={(410.09,113.92)}, rotate = 270.44] [fill={rgb, 255:red, 0; green, 0; blue, 0 }  ][line width=0.08]  [draw opacity=0] (8.75,-4.2) -- (0,0) -- (8.75,4.2) -- (5.81,0) -- cycle    ;
\draw [line width=1.5]    (196.22,172.04) -- (172.07,149.28) ;
\draw [shift={(169.16,146.54)}, rotate = 43.29] [fill={rgb, 255:red, 0; green, 0; blue, 0 }  ][line width=0.08]  [draw opacity=0] (8.75,-4.2) -- (0,0) -- (8.75,4.2) -- (5.81,0) -- cycle    ;
\draw [line width=1.5]    (298.42,60.42) -- (274.27,37.66) ;
\draw [shift={(271.35,34.92)}, rotate = 43.29] [fill={rgb, 255:red, 0; green, 0; blue, 0 }  ][line width=0.08]  [draw opacity=0] (8.75,-4.2) -- (0,0) -- (8.75,4.2) -- (5.81,0) -- cycle    ;
\draw   (192.75,92.55) .. controls (192.54,89.99) and (194.46,87.75) .. (197.02,87.55) .. controls (199.58,87.35) and (201.82,89.26) .. (202.02,91.83) .. controls (202.22,94.39) and (200.31,96.63) .. (197.75,96.83) .. controls (195.18,97.03) and (192.95,95.11) .. (192.75,92.55) -- cycle ;
\draw    (151.07,49) -- (193.41,88.83) ;
\draw [shift={(149.61,47.63)}, rotate = 43.25] [fill={rgb, 255:red, 0; green, 0; blue, 0 }  ][line width=0.08]  [draw opacity=0] (12,-3) -- (0,0) -- (12,3) -- cycle    ;
\draw    (239.68,45.24) -- (200.04,87.77) ;
\draw [shift={(241.05,43.78)}, rotate = 132.99] [fill={rgb, 255:red, 0; green, 0; blue, 0 }  ][line width=0.08]  [draw opacity=0] (12,-3) -- (0,0) -- (12,3) -- cycle    ;
\draw  [fill={rgb, 255:red, 0; green, 0; blue, 0 }  ,fill opacity=1 ] (196.51,93.1) .. controls (196,92.62) and (195.99,91.82) .. (196.47,91.31) .. controls (196.95,90.81) and (197.75,90.79) .. (198.26,91.27) .. controls (198.77,91.76) and (198.78,92.56) .. (198.3,93.07) .. controls (197.82,93.57) and (197.01,93.59) .. (196.51,93.1) -- cycle ;
\draw  [draw opacity=0] (295.51,113.84) .. controls (295.56,113.84) and (295.62,113.84) .. (295.67,113.84) .. controls (301.93,113.99) and (306.85,120.97) .. (306.66,129.45) .. controls (306.57,133.33) and (305.43,136.85) .. (303.61,139.52) -- (295.32,129.19) -- cycle ; \draw    (295.67,113.84) .. controls (301.93,113.99) and (306.85,120.97) .. (306.66,129.45) .. controls (306.59,132.56) and (305.85,135.43) .. (304.62,137.82) ; \draw [shift={(303.61,139.52)}, rotate = 291.91] [fill={rgb, 255:red, 0; green, 0; blue, 0 }  ][line width=0.08]  [draw opacity=0] (7.2,-1.8) -- (0,0) -- (7.2,1.8) -- cycle    ; \draw [shift={(295.51,113.84)}, rotate = 24.68] [fill={rgb, 255:red, 0; green, 0; blue, 0 }  ][line width=0.08]  [draw opacity=0] (7.2,-1.8) -- (0,0) -- (7.2,1.8) -- cycle    ;

\draw (413.67,79.33) node [anchor=north west][inner sep=0.75pt]   [align=left] {{\small $\vect g$}};
\draw (453,190) node [anchor=north west][inner sep=0.75pt]   [align=left] {{\small $\hat \imath _{\rm{A}}$}};
\draw (375,120) node [anchor=north west][inner sep=0.75pt]   [align=left] {{\small $\hat \jmath _{\rm{A}}$}};
\draw (234,31) node [anchor=north west][inner sep=0.75pt]  [rotate=-313.75] [align=left] {{\small $\hat \imath _{\rm B}$}};
\draw (129.9,40) node [anchor=north west][inner sep=0.75pt]  [rotate=-313.75] [align=left] {{\small $\hat \jmath _{\rm B}$}};
\draw (287,20.33) node [anchor=north west][inner sep=0.75pt]   [align=left] {{\small ${\vect f}_2$}};
\draw (184,135.33) node [anchor=north west][inner sep=0.75pt]   [align=left] {{\small ${\vect f}_1$}};
\draw (260.67,125) node [anchor=north west][inner sep=0.75pt]   [align=left] {{\small c}};
\draw (327.33,200) node [anchor=north west][inner sep=0.75pt]   [align=left] {{\small \textit{w}}};
\draw (308.33,115) node [anchor=north west][inner sep=0.75pt]   [align=left] {{\small $\phi$}};
\end{tikzpicture}
}
\caption{Bicopter configuration considered in this paper. The bicopter is constrained to the {$\hat  \imath _\rmA- \hat \jmath _\rmA$} plane and rotates about the {$\hat k_\rmA$} axis of the inertial frame {$\rm F_A.$}
Note that $\vect r_{\rmc/w} = r_1 \hat \imath_\rmA + r_2 \hat \jmath_\rmA.$}
    \label{fig:Bicopter}
\end{figure}

We consider the bicopter system shown in Figure \ref{fig:Bicopter}.
The derivation of the equations of motion and their dynamic extension is described in more detail in \cite{portella2024adaptive}.

The dynamically extended equations of motion are
\begin{align}
    \dot x_1 &= x_2,
    \label{eq:x1dot}
    \\
    \dot x_2 &= f_2(x_1, x_2) + g_2(x_3) \theta_1 ,
    \label{eq:x2dot}
    \\
    \dot x_3 &= x_4,
    \label{eq:x3dot}
    \\
    \dot x_4
        &=
            g_4(\theta)
            u \label{eq:x4dot},
\end{align}
where
\begin{align}
    x_1 \isdef \matl
                 r_1\\
                 r_2
               \matr, 
    x_2 \isdef \matl
                 \dot r_1\\
                 \dot r_2
               \matr,
    x_3 \isdef \matl
                 F\\
                 \phi
               \matr, 
    x_4 \isdef \matl
                 \dot F\\
                 \dot \phi
               \matr,
\end{align}
and 
\begin{align}
    u
        \isdef 
            \matl 
                \ddot F \\
                M
            \matr,
    \theta
        \isdef
            \matl
                m^{-1}\\
                J^{-1}
            \matr.
\end{align}
$r_1$ and $r_2$ are the horizontal and vertical positions of the center of mass of the bicopter, 
$\dot r_1$ and $\dot r_2$ are the corresponding velocities, 
$\phi$ is the roll angle of the bicopter, 
$\dot \phi$ is the angular velocity of the bicopter, 
$F$ and $M$ are the total force and the total moment applied to the bicopter, respectively,  
and 
$m$ and $J$ are the mass and the moment of inertia of the bicopter, respectively. 
The functions 
\begin{align}
    f_2(x_1, x_2) 
        &\isdef 
            \matl
                0\\
                -g
            \matr,
            \\ 
    g_2(x_3) 
        &\isdef
            \matl
                -\sin(x_{3,2})x_{3,1}\\
                \cos(x_{3,2}) x_{3,1}
            \matr, \\
    g_4(\theta)
        &\isdef 
            \matl
                0 & \theta_2\\
                1 & 0
            \matr.            
\end{align}
Note that the system \eqref{eq:x1dot}-\eqref{eq:x4dot} is in pure feedback form since $x_3$ appears non-affinely in \eqref{eq:x2dot}. 
Finally, since $g_4(\theta)$ is invertible, the backstepping framework can be applied to design a stabilizing controller.


\section{Adaptive Backstepping Control}
\label{sec:adaptive_backstepping}
In this section, we develop an adaptive controller based on the backstepping framework to stabilize the bicopter system. 
Unlike the controller designed in \cite{portella2024adaptive}, the proposed controller does not require inversion of the parameter estimates. 


Let $\xi_1 \isdef \matl r_{\rmd1} & r_{\rmd2} \matr^\rmT \in \BBR^2$ denote the desired trajectory and define the tracking error $e_1 \isdef x_1 - \xi_1.$

\subsection{$e_1$ stabilization}
Consider the function
\begin{align}
    V_1
        &\isdef
            \dfrac{1}{2}e_1^{\rm T}e_1.
            \label{eq:V1}
\end{align}
Differentiating \eqref{eq:V1} and using \eqref{eq:x1dot} yields
$\dot V_1
        =
            e_1^\rmT \dot x_1
        =
            e_1^\rmT x_2 .$
Note that if $x_2 = -k_1 e_1,$ where $k_1>0,$ then $\dot V_1 < 0.$ 
However, $x_2$ is not the control signal and thus cannot be chosen arbitrarily.
Instead, following the backstepping process, we design a control law that yields the desired $x_2$ response, as follows.

Define
$
    e_2
        \isdef
            x_2 - \xi_2,
$
where
\begin{align}
    \xi_2 
        \isdef 
            -k_1e_1,
    \label{eq:xi_2}
\end{align} 
and $k_1>0$. Then, differentiating \eqref{eq:V1} yields
\begin{align}
    \dot{V}_1 
        &= 
            -k_1e_1^{\rm T}e_1 + e_2^{\rm T}e_1.
    \label{eq:V_1_dot}
\end{align}

\subsection{$e_2$ stabilization}
Define $p_1 \isdef \theta_1^{-1}$ and consider the function
\begin{align}
    V_2
        &\isdef
            V_1 + \dfrac{1}{2}e_2^{\rm T}e_2 + \dfrac{1}{2}\gamma_1^{-1}|\theta_1|\left(\hat{p}_1 - p_1\right)^2,
    \label{eq:V_2}
\end{align}
where $\gamma_1>0$ and $\hat{p}_1$ is an estimate of $p_1$. 
Differentiating \eqref{eq:V_2} and using \eqref{eq:V_1_dot} yields 
$
    \dot{V}_2 
        = 
            -k_1e_1^{\rm T}e_1 
            + e_2^{\rm T}\left(e_1 + f_2 + g_2(x_3)\theta_1 + k_1x_2\right) 
            + \gamma_1^{-1}|\theta_1|\left(\hat{p}_1 - p_1\right)\dot{\hat{p}}_1.
$ 
However, $g(x_3)$ is not the control signal and thus can not be chosen arbitrarily.
Instead, following the backstepping process, we design a control law that yields the desired $g(x_3)$ response, as follows.

Define
$
    e_3
        \isdef
            g(x_3) - \xi_3,
$
where
\begin{align}
    \xi_3 
        \isdef
            -\hat{p}_1\left(e_1 + f_2 + k_1x_2\right) 
            - k_2\sigma(\theta_1)e_2,
        \label{eq:xi_3}
\end{align}
$\sigma(\theta_1)$ is the sign of $\theta_1,$
and $k_2>0$. 
Consider the parameter update law 
\begin{align}
    \dot{\hat{p}}_1 
        = 
            \gamma_1\sigma(\theta_1)e_2^{\rm T}\left(e_1 + f_2 + k_1x_2\right).
        \label{eq:p1_hat_dot}
\end{align}
Then, differentiating \eqref{eq:V_2} yields
\begin{align}
    \dot{V}_2
        &=
            -k_1e_1^{\rm T}e_1 - k_2|\theta_1|e_2^{\rm T}e_2 + e_3^{\rm T}e_2\theta_1.
    \label{eq:V_2_dot}
\end{align}

\subsection{$e_3$ stabilization} 
Consider the function
\begin{align}
    V_3
        &\isdef
            V_2 + \dfrac{1}{2}e_3^{\rm T}e_3.
    \label{eq:V_3}
\end{align}
Differentiating \eqref{eq:V_3} and using \eqref{eq:V_2_dot} yields 
$
    \dot{V}_2 
        = 
            -k_1e_1^{\rm T}e_1 
            - k_2|\theta_1|e_2^{\rm T}e_2 
            + e_3^{\rm T}(\theta_1e_2 
            + \dfrac{\partial g(x_3)}{\partial x_3}x_4 
            + \dot{\hat{p}}_1\left(e_1 + f_2 + k_1x_2\right) 
            + \hat{p}_1\left(x_2 + k_1f_2\right) 
            + k_2\sigma(\theta_1)\left(f_2 + k_1x_2\right) 
            + \left(k_1\hat{p}_1 + k_2\sigma(\theta_1)\right)g_2(x_3)\theta_1).
$
However, $x_4$ is not the control signal and thus can not be chosen arbitrarily.
Instead, following the backstepping process, we design a control law that yields the desired $x_4$ response, as follows.

Define
$
    e_4
        \isdef
            x_4 - \xi_4,
$
where
\begin{align}
    \xi_4 
        &\isdef 
            \SG_2^{-1}
            \Big[\dot{\hat{p}}_1\left(e_1 + f_2 + k_1x_2\right) 
            + \hat{p}_1\left(x_2 + k_1f_2\right)
            \nn\\ &\quad
            + k_2\sigma(\theta_1)\left(f_2 + k_1x_2\right)
            \nn\\ &\quad
            + \left(e_2 + \left(k_1\hat{p}_1 + k_2\sigma(\theta_1)\right)g_2(x_3)\right)\hat{\theta}_1 
            + k_3e_3\Big],
    \label{eq:xi_4}
\end{align}
$k_3>0$ and  

\begin{align}
    \SG_2 \isdef \frac{\partial g_2}{\partial x_3}
            =
            \matl
                -\cos(x_{3,2})x_{3,1} & -\sin(x_{3,2})\\
                -\sin(x_{3,2})x_{3,1} & \cos(x_{3,2})
            \matr.
\end{align}

Then, differentiating \eqref{eq:V_3} yields
\begin{align}
    \dot{V}_3
        &=
            - k_1e_1^{\rm T}e_1 
            - k_2|\theta_1|e_2^{\rm T}e_2
            - k_3e_3^{\rm T}e_3
            \nn\\ &\quad
            + e_4^{\rm T}\SG_2^{\rm T}e_3
            + (\theta_1 - \hat{\theta}_1)e_3^{\rm T}(e_2
            \nn \\ &\quad
            + \left(k_1\hat{p}_1 + k_2\sigma(\theta_1)\right)g_2(x_3)).
    \label{eq:V3_dot}
\end{align}

\subsection{$e_4$ stabilization}
Define $p_2 \isdef \theta_2^{-1}$ and consider the function
\begin{align}
    V_4
        &\isdef
            V_3
            + \dfrac{1}{2}e_4^{\rm T}e_4
            + \dfrac{1}{2}\gamma_2^{-1}|\theta_2|\left(\hat{p}_2 - p_2\right)^2
            \nn\\ &\quad
            + \dfrac{1}{2}\alpha_1^{-1}
            (\theta_1 - \hat{\theta}_1)^2
            + \dfrac{1}{2}\alpha_2^{-1}
            (\theta_1 - \hat{\vartheta}_1)^2,
    \label{eq:v_4}
\end{align}
where $\gamma_2>0$, $\alpha_1>0$, $\alpha_2>0,$ 
$\hat{p}_2$ is an estimate of $p_2$, and  
$\hat{\theta}_1$ and $\hat \vartheta_1$ are estimates of $\theta_1$.

Differentiating \eqref{eq:v_4} 
yields
\begin{align}
    \dot{V}_4
        &=
            - k_1e_1^{\rm T}e_1 
            - k_2|\theta_1|e_2^{\rm T}e_2
            - k_3e_3^{\rm T}e_3
            + 
            \nn \\ &\quad 
            (\theta_1 - \hat{\theta}_1)
            e_3^{\rm T}(e_2
            + \left(k_1\hat{p}_1 + k_2\sigma(\theta_1)\right)g_2(x_3) 
            )
            -
            \nn \\ &\quad 
            (\theta_1 - \hat{\theta}_1)
            e_3^{\rm T}
            \alpha_1^{-1}\dot{\hat{\theta}}_1
            +
            \nn \\ &\quad 
            e_4^{\rm T}
            (g_4(\theta) u + \Sigma + \Psi\theta_1)
            +
            \nn \\ &\quad 
            \gamma_2^{-1}|\theta_2|\left(\hat{p}_2 - p_2\right)\dot{\hat{p}}_2
            - \alpha_2^{-1}
            (\theta_1 - \hat{\vartheta}_1)
            \dot{\hat{\vartheta}}_1, 
    \label{eq:V_4_dot}
\end{align}
where 
\begin{align}
    \Sigma
        &\isdef
            d_t\text{inv}\SG_2
            \Bigg[\dot{\hat{p}}_1\left(e_1 + k_1x_2 + f_2\right)
            \nn\\ &\quad
            + \hat{p}_1\left(x_2 + k_1f_2\right)
            + k_2\sigma(\theta_1)\left(f_2 + k_1x_2\right)
            \nn\\ &\quad
            + \hat{\theta}_1\left(e_2 + k_1\hat{p}_1g_2(x_3) + k_2\sigma(\theta_1)g_2(x_3)\right) + k_3e_3\Bigg]
            \nn\\ &\quad
            + \SG_2^{-1}\Bigg[\gamma_1\sigma(\theta_1)\Big[\left(f_2^{\rm T} + k_1x_2^{\rm T}\right)\left(e_1 + k_1x_2 + f_2\right)
            \nn\\ &\quad
            + e_2^{\rm T}\left(x_2 + k_1f_2\right)\Big]\left(e_1 + k_1x_2 + f_2\right)
            \nn\\ &\quad
            + 2\dot{\hat{p}}_1\left(x_2 + k_1f_2\right)
            \nn\\ &\quad
            + 2\hat{p}_1f_2 + k_2\sigma(\theta_1)k_1f_2 
            + \dot{\hat{\theta}}_1\Big(e_2 + k_1\hat{p}_1g_2(x_3) 
            \nn\\ &\quad
            + k_2\sigma(\theta_1)g_2(x_3)\Big) + \hat{\theta}_1\Big(f_2 + k_1x_2 + k_1\hat{p}_1g_2(x_3) 
            \nn\\ &\quad
            + k_1\hat{p}_1\SG_2x_4 + k_2\sigma(\theta_1)\SG_2x_4\Big)
            \nn\\ &\quad
            + k_3\Big(\SG_2x_4 + \dot{\hat{p}}_1\left(e_1 + k_1x_2 + f_2\right) + \hat{p}_1\left(x_2 + k_1f_2\right)
            \nn\\ &\quad
            + k_2\sigma(\theta_1)\left(f_2 + k_1x_2\right)\Big)\Bigg]
            \in \BBR^2,
        \label{eq:Sigma}
%
%
    \\
    \Psi
        &\isdef
            \SG_2^{-1}
            \Bigg[\gamma_1\sigma(\theta_1)(g_2(x_3)^{\rm T}\Phi)\Phi
            + \dot{\hat{p}}_1k_1g_2(x_3) 
            \nn\\ &\quad
            + \hat{p}_1g_2(x_3)
            + k_1k_2\sigma(\theta_1)g_2(x_3)
            + \hat{\theta}_1g_2(x_3) 
            \nn\\ &\quad
            + k_1k_3\hat{p}_1g_2(x_3) 
            + k_2k_3\sigma(\theta_1)g_2(x_3)\Bigg]
            \in \BBR^2
        \label{eq:psi}, 
    \\
    \Phi 
        &\isdef
            e_1 + k_1x_2 + f_2
            \in \BBR^2.
\end{align}

Finally, consider the control law
\begin{align}
    u
        &=
            -\left[
                \matl
                    0 & 1\\
                    \hat{p}_2 & 0
                \matr
                (\Sigma + \Psi\hat{\vartheta}_1)
                + k_4
                \matl
                    0 & 1\\
                    \sigma(\theta_2) & 0
                \matr
                e_4
            \right],
        \label{eq:u}
\end{align}
where $k_4>0$ and the parameters update laws 
\begin{align}
    \dot{\hat{p}}_2
        &=
            \gamma_2\sigma(\theta_2)
            \matl
                1 & 0
            \matr
            e_4
            \matl
                1 & 0
            \matr
            (\Sigma + \Psi\hat{\vartheta}_1),
        \label{eq:p2_hat_dot}
    \\
    \dot{\hat{\theta}}_1
        &=
            \alpha_1e_3^{\rm T}\left(e_2 + k_1\hat{p}_1g_2(x_3) + k_2\sigma(\theta_1)g_2(x_3)\right),
        \label{eq:theta_1_hat_dot}
    \\
    \dot{\hat{\vartheta}}_1
        &=
            \alpha_2\Sigma^{\rm T}e_4.
        \label{eq:vartheta_1_hat_dot}
\end{align}

Substituting the control law and the parameter update laws above in \eqref{eq:V_4_dot} yields
\begin{align}
    \dot{V}_4
        &=
            - k_1e_1^{\rm T}e_1 
            - k_2|\theta_1|e_2^{\rm T}e_2
            \nn\\ &\quad
            - k_3e_3^{\rm T}e_3
            - k_4e_4^{\rm T}
            \matl
                |\theta_2| & 0\\
                 0 & 1
            \matr
            e_4
        \leq
            0.
        \label{eq:V_final_dot}
\end{align}

\subsection{Stability Analysis}

The following theorem uses the Barbashin-Krasovskii-LaSalle invariance principle to establish the asymptotic stability of the origin in the closed-loop system \eqref{eq:x1dot}, \eqref{eq:x2dot}, \eqref{eq:x3dot}, \eqref{eq:x4dot}, \eqref{eq:u}.

\begin{theorem}
    Consider the system \eqref{eq:x1dot}, \eqref{eq:x2dot}, \eqref{eq:x3dot}, \eqref{eq:x4dot}, the controller \eqref{eq:u}, and the parameter update laws \eqref{eq:p1_hat_dot}, \eqref{eq:p2_hat_dot}, \eqref{eq:theta_1_hat_dot}, \eqref{eq:vartheta_1_hat_dot}, 
    where 
    $\gamma_1, \gamma_2, \alpha_1, \alpha_2, \alpha_3, k_1,  k_2, k_3, k_4 > 0. $
    Then, the equilibrium point $(e_1,e_2,e_3,e_4) = (0,0,0,0)$ is asymptotically stable.
\end{theorem}

\begin{proof}
Consider the function \eqref{eq:v_4}.
Note that, for all $e_1, e_2, e_3, e_4 \in \BBR^2,$ it follows from \eqref{eq:V_4_dot} that $\dot V_4 \leq 0.$
%
Let $\ell >0$ and define
\begin{align}
    \Omega
        &\isdef
            \left\{
                e_1,e_2,e_3,e_4 \in \mathbb{R}^{2},
                \hat{p}_1,\hat{p}_2,\hat{\theta}_1,\hat{\vartheta}_1 \in \mathbb{R}: V_4 \leq \ell
            \right\}. \nn
\end{align}
Note that $\Omega $ is a compact set.
Furthermore, since $\dot V_4 \leq 0 $ and $V_4 > 0,$ $\Omega$ is a {invariant} set.

Next, define $E \subset \Omega$ such that $e_1 =  e_2 = e_3 = e_4 = 0.$
Since $\dot{V_4} = 0$ if and only if at $e_1 = e_2 = e_3 = e_4 = 0,$ it follows that $E$ is the largest positively invariant set such that $\dot V_4 = 0.$ 
It thus follows from the Barbashin-Krasovskii-LaSalles's invariance principle \cite[p.~129]{khalil2002nonlinear} that $(e_1,e_2,e_3,e_4) = (0,0,0,0)$ is asymptotically stable.
\end{proof}

\begin{remark}
    Note that the control \eqref{eq:u} requires the existence of $\Sigma$ and $\Psi$, which implies that $\SG_2$ must be non-singular. If $F \neq 0,$ which is reasonably expected during the system's operation, $\SG_2$ is nonsingular. 
\end{remark}

\section{Numerical Simulation}
\label{sec:simulations}

\begin{figure}[h]
    \centering
    \vspace{-0.5em}
        \resizebox{\columnwidth}{!}{%
            \begin{tikzpicture}[>={stealth'}, line width = 0.25mm]
                \node [smallblock, fill=green!20,  rounded corners, minimum height = 5cm, minimum width = 6cm] (Parameter_estimator) {};

                \node [input,left = 2.3cm of Parameter_estimator] (Input_parameter_estimator) {};

                \node [output, right = 2.3 cm of Parameter_estimator] (Output_parameter_estimator) {};

                \draw[-] (Input_parameter_estimator) --node[above]
                {$(x_1,x_2,x_3,x_4)$} (Parameter_estimator.west);

                \draw[->]  (Parameter_estimator.east) --node[above]
                {$(\hat{\theta}_1,\hat{\vartheta}_1,\hat{p}_1,\hat{p}_2)$} (Output_parameter_estimator);

                \node [smallblock, rounded corners, right = 3.0cm of Parameter_estimator.west, minimum height = 0.6cm, minimum width = 0.7cm, yshift = 1.8cm](p_hat_1){$\eqref{eq:p1_hat_dot}$};
                
                \node [smallblock, rounded corners, below = 0.5cm of p_hat_1, minimum height = 0.6cm, minimum width = 0.7cm](p_hat_2){$\eqref{eq:p2_hat_dot}$};

                \node [smallblock, rounded corners, below = 0.5cm of p_hat_2, minimum height = 0.6cm, minimum width = 0.7cm](theta_hat_1){$\eqref{eq:theta_1_hat_dot}$};

                \node [smallblock, rounded corners, below = 0.5cm of theta_hat_1, minimum height = 0.6cm, minimum width = 0.7cm](vartheta_hat_1){$\eqref{eq:vartheta_1_hat_dot}$};


                \node [smallblock, rounded corners, right = 5.0 of Parameter_estimator.west, minimum height = 0.6cm, minimum width = 0.7cm](integral){$\dfrac{1}{s}$};


                \draw [->] (Parameter_estimator.west) -| ($(Parameter_estimator.west) + (0.5cm, 0)$) |-node[near end, above]{$(x_1,x_2)$} (p_hat_1.west);

                \draw [->] (Parameter_estimator.west) -| ($(Parameter_estimator.west) + (0.5cm, 0)$) |-node[near end, above]{$(x_1,x_2,x_3,x_4)$} (p_hat_2.west);

                \draw [->] (Parameter_estimator.west) -| ($(Parameter_estimator.west) + (0.5cm, 0)$) |-node[near end, above]{$(x_1,x_2,x_3)$} (theta_hat_1.west);

                \draw [->] (Parameter_estimator.west) -| ($(Parameter_estimator.west) + (0.5cm, 0)$) |-node[near end, above]{$(x_1,x_2,x_3,x_4)$} (vartheta_hat_1.west);


                \draw [-] (p_hat_1.east) -| node[near start, above]{$\dot{\hat{p}}_1$} ($(Parameter_estimator.east) + (-1.2cm, 0)$) -- (integral.west);

                \draw [-] (p_hat_2.east) -| node[near start, above]{$\dot{\hat{p}}_2$} ($(Parameter_estimator.east) + (-1.2cm, 0)$) -- (integral.west);

                \draw [-] (theta_hat_1.east) -| node[near start, above]{$\dot{\hat{\theta}}_1$} ($(Parameter_estimator.east) + (-1.2cm, 0)$) -- (integral.west);

                \draw [-] (vartheta_hat_1.east) -| node[near start, above]{$\dot{\hat{\vartheta}}_1$} ($(Parameter_estimator.east) + (-1.2cm, 0)$) -- (integral.west);

                \draw [-] (integral.east) -- (Parameter_estimator.east);
            \end{tikzpicture}
        }
    \caption{
    Parameter estimation architecture. 
    }
    \label{Parameter_estimators_block_diag}
\end{figure}
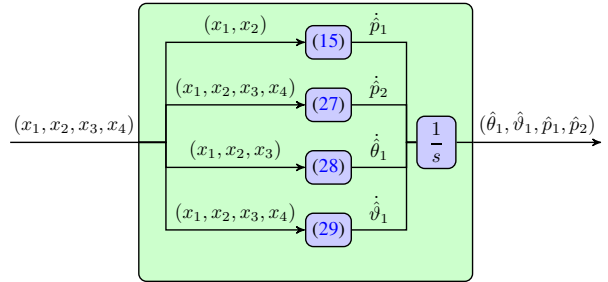

In this section, we apply the adaptive backstepping controller developed in the previous section to the trajectory tracking problem.
In particular, the adaptive controller is used to follow an elliptical and a second-order Hilbert curve-based trajectory. 
Note that the adaptive controller is given by \eqref{eq:u} with the parameter update laws \eqref{eq:p1_hat_dot}, \eqref{eq:p2_hat_dot}, \eqref{eq:theta_1_hat_dot}, and \eqref{eq:vartheta_1_hat_dot}.
The control architecture is shown in Figure \ref{Adaptive_backstepping_blk_diag} and the estimation architecture is shown in Figure \ref{Parameter_estimators_block_diag}.

\begin{figure}[h]
\centering
\vspace{-1em}
    \resizebox{0.9\columnwidth}{!}{%
    \begin{tikzpicture}[>={stealth'}, line width = 0.25mm]

    \node [input, name=ref]{};
    \node [smallblock, rounded corners, right = 0.5cm of ref , minimum height = 0.6cm, minimum width = 0.7cm] (controller) {$\begin{array}{c} {\rm Adaptive} \\ {\rm Controller} \\ \eqref{eq:u}\end{array}$};
    \node [smallblock, rounded corners, right = 0.75cm of controller, minimum height = 0.6cm , minimum width = 0.5cm] (system) {$\begin{array}{c} {\rm Bicopter} \\ {\rm System}\end{array}$};

    \node [smallblock, fill=green!20, rounded corners, below = 1.5 cm of controller, minimum height = 0.6cm , minimum width = 0.5cm] (estimator) {$\begin{array}{c} {\rm Parameter} \\ {\rm Estimator}\end{array}$};
    
    \node [output, right = 1.0cm of system] (output) {};
    \node [input, below = 2.25cm of system] (midpoint) {};
    \node [input, below = 0.75cm of controller] (midpoint2) {};
    \node [input, left = 1.0cm of controller] (reference) {};
    
    \draw [->] (controller) -- node [above] {$u (t)$} (system);
    
    
    \draw [->] (system) -- node [name=y, near end]{} node [very near end, above] {}(output);
    
    \draw [->] (reference.west) -- node[xshift=-1em, yshift=2em]{$\begin{array}{c} {\rm Desired}\\ {\rm Trajectory} \end{array}$} (controller.west);
    \draw [-] (y.west) |- (midpoint);
    \draw [->] (midpoint) -- node[near start,above] {$(x_1,x_2,x_3,x_4)$} (estimator.east);
    \draw [->] ([xshift = -0.75cm]estimator.north) -- node [near start, right] {$(\hat{\theta}_1,\hat{\vartheta}_1,\hat{p}_1,\hat{p}_2)$} ([xshift = -0.75cm]controller.south);
    \draw [-] ([xshift = -0.0cm]y.west) |- (midpoint2);
    \draw [->] (midpoint2) -- node [right] {$(x_1,x_2,x_3,x_4)$} (controller.south);
    
    \end{tikzpicture}
    }  
    \caption{
        Singularity-free adaptive control architecture. 
    }
    \vspace{-1em}
    \label{Adaptive_backstepping_blk_diag}
\end{figure}
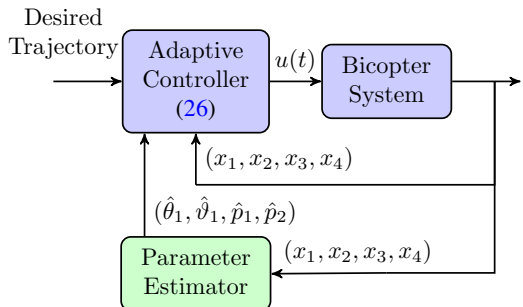

To simulate the bicopter, we assume that the mass of the bicopter $m$ is $1 \ \rm kg$ and its moment of inertia $J$ is $0.2$ $\rm kg  m^2.$
In the controller, 
we set 
$k_1 = 1$, $k_2 = 5$, $k_3 = 10$, and $ k_4 = 20. $
In the adaptation laws presented in Figure \ref{Parameter_estimators_block_diag}, 
we set the estimator gains 
$\gamma_1 = \gamma_2 =  \alpha_1 = \alpha_2 = 0.1$,  
and the initial estimates 
$\hat{\theta}_1(0) = \hat{\vartheta}_1(0) = \hat{p}_1(0) = \hat{p}_2(0) = 0.$

\subsection{Elliptical Trajectory}
The bicopter is commanded to follow an elliptical trajectory given by 
{
\begin{align}
    r_{\rmd1}(t) 
        &=
            5 \cos(\psi)-5 \cos(\psi) \cos (\omega t) - 3 \sin(\psi)\sin(\omega t), 
    \\
    r_{\rmd2}(t) 
        &=
            5\sin(\psi) - 5 \sin(\psi) \cos (\omega t) + 3 \cos(\psi)\sin(\omega t), 
\end{align}
where $\psi=45~\rm{deg}$ and $\omega = 0.1 \ \rm rad/s^{-1}. $
Figure \ref{fig:Elliptical trajectory} shows the trajectory-tracking response of the bicopter, where the desired trajectory is shown in black dashes, and the output trajectory response is shown in blue.
Figure \ref{fig:States Elliptical} shows the position $r_1, r_2$ response, the roll angle $\phi$ response of the bicopter and the control $u$ generated by the adaptive backstepping controller \eqref{eq:u}. 
%
%

\begin{figure}[!ht]
    \centering
    \includegraphics[width=\linewidth]{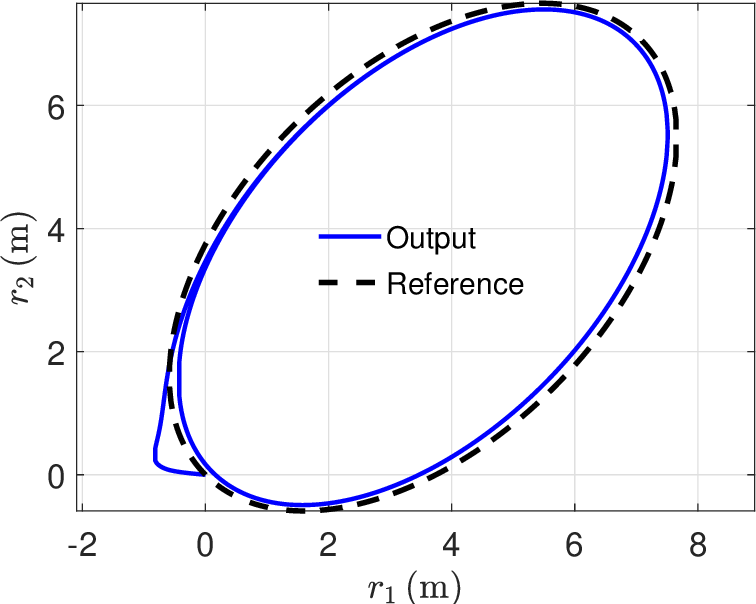}
    \caption{\textbf{Elliptical trajectory}. Tracking response of the bicopter with the adaptive backstepping controller. 
    Note that the output trajectory is shown in solid blue, and the reference trajectory is in dashed black.}
    \label{fig:Elliptical trajectory}
\end{figure}

\begin{figure}[!ht]
    \centering
    \includegraphics[width=\linewidth]{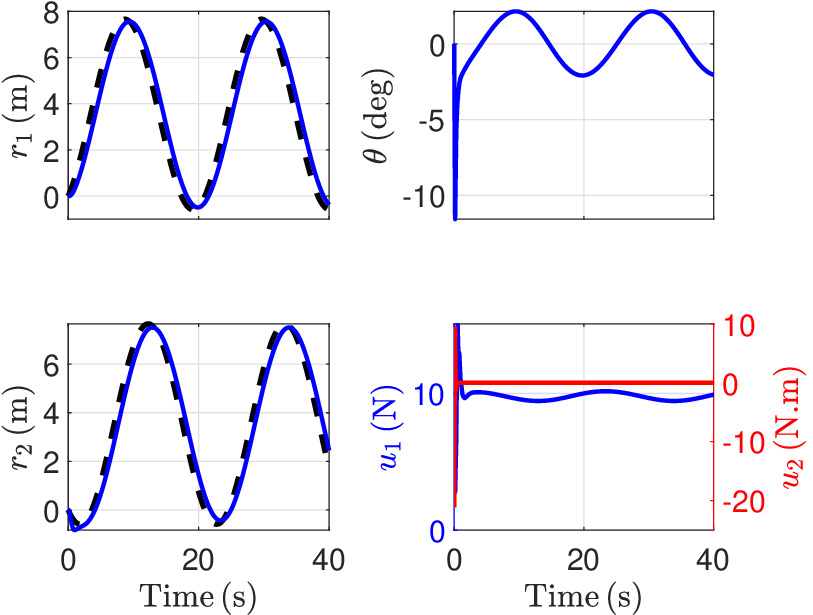}
    \caption{\textbf{Elliptical trajectory}. Position $(r_1, r_2)$, roll angle $\phi$ response and the input signal applied to the bicopter obtained with adaptive backstepping controller \eqref{eq:u}. Note that the references $r_{\rmd1},r_{\rmd2}$ is shown in black dashed lines}
    \label{fig:States Elliptical}
\end{figure}


Figure \ref{fig:Parameter_estimations_Elliptical} shows the estimates $\hat \theta_1,\hat \vartheta_1,$ and $ \hat{p}_1$ of $\theta_1^{-1}$ and the estimate $\hat{p}_2$ of $\theta_2^{-1}.$
Note that the parameter estimates do not converge to their actual values.
However, the non-convergence of the estimates is not due to persistency-related issues. 
In the adaptive controller design, since the parameter adaptation laws are chosen to cancel undesirable factors and not to estimate the parameters, the estimates do not necessarily need to converge.

\begin{figure}[!ht]
    \centering
    \includegraphics[width=\linewidth]{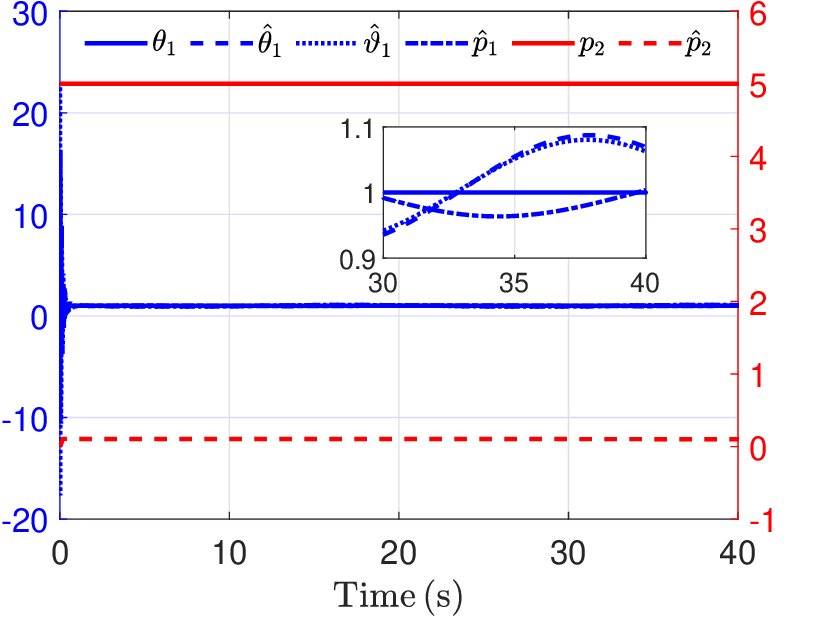}
    \caption{\textbf{Elliptical trajectory}. 
    Estimates of $\theta_1$, $p_1 = \theta_1^{-1}$, and $p_2 = \theta_2^{-1}$ obtained with adaption laws 
    \eqref{eq:theta_1_hat_dot}, 
    \eqref{eq:vartheta_1_hat_dot},
    \eqref{eq:p1_hat_dot}, and
    \eqref{eq:p2_hat_dot}.}
    \label{fig:Parameter_estimations_Elliptical}
\end{figure}

Next, to improve the tracking performance, we increase the estimator gains by a factor of 10; that is, the gains are $\gamma_1 = \gamma_2 =  \alpha_1 = \alpha_2 = 1.$  
Figures \ref{fig:elliptical_trajectory_faster_estimation}, 
\ref{fig:elliptical_trajectory_states_faster_estimation}, 
and 
\ref{fig:elliptical_parameter_estimates_faster_estimation} show the trajectory tracking response, input signals, and the parameter estimates with the larger estimator gains. 
Note that, although the trajectory tracking error improves, the required control effort also increases. 

\begin{figure}[!ht]
    \centering
    \includegraphics[width=\linewidth]{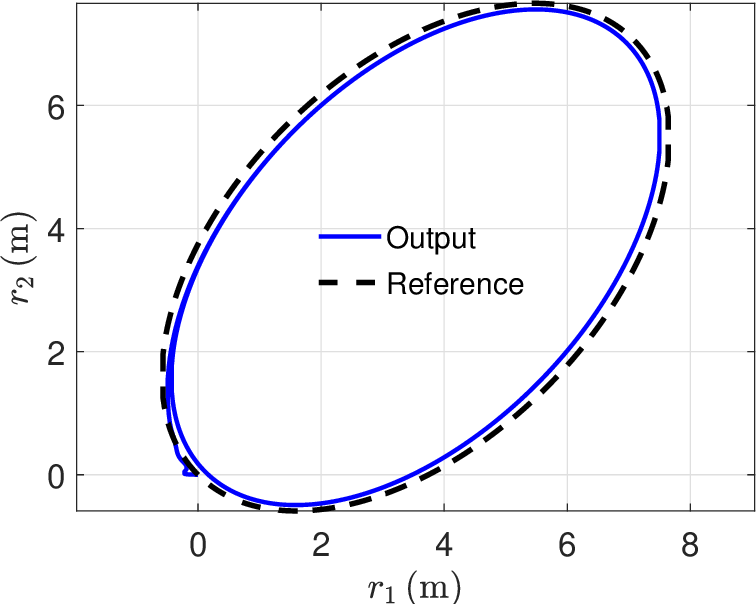}
    \caption{\textbf{Elliptical trajectory}. Tracking response of the bicopter with faster parameter estimation. 
    Note that the output trajectory is shown in solid blue, and the reference trajectory is in dashed black.}
    \label{fig:elliptical_trajectory_faster_estimation}
\end{figure}

\begin{figure}[!ht]
    \centering
    \includegraphics[width=\linewidth]{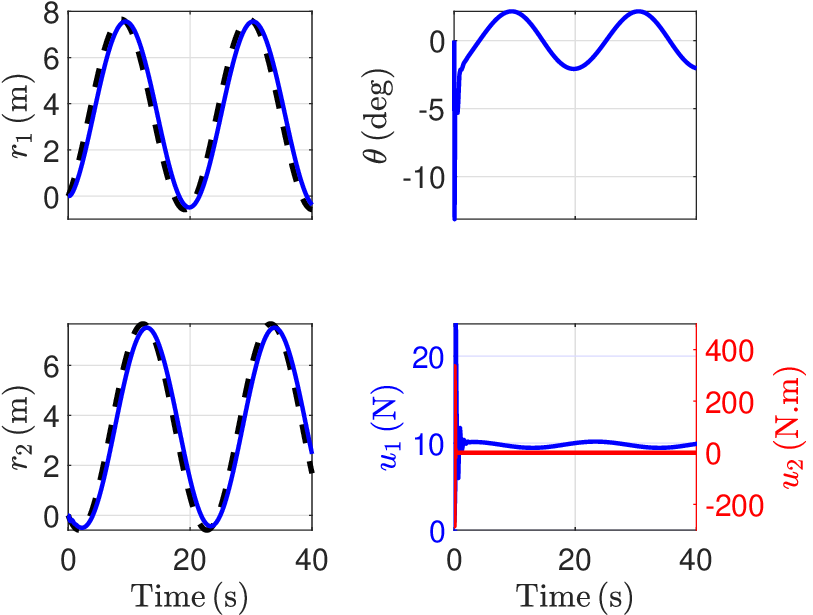}
    \caption{\textbf{Elliptical trajectory}. Position $(r_1, r_2)$ and roll angle $\phi$ response of the bicopter obtained with faster parameter estimation. Note that the references $r_{\rmd1},r_{\rmd2}$ is shown in black dashed lines}
    \label{fig:elliptical_trajectory_states_faster_estimation}
\end{figure}

\begin{figure}[!ht]
    \centering
    \includegraphics[width=\linewidth]{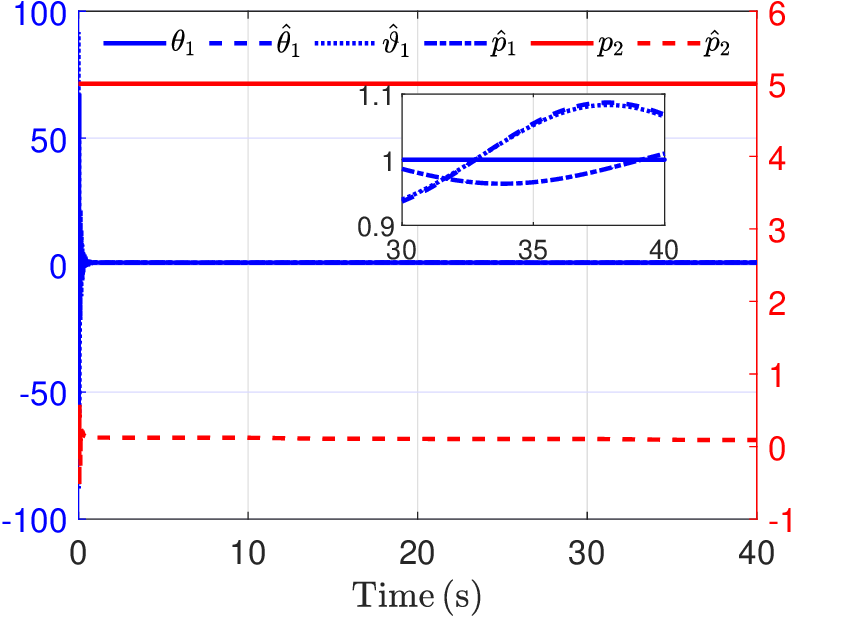}
    \caption{\textbf{Elliptical trajectory}. Estimates of $\theta_1$, $p_1 = \theta_1^{-1}$, and $p_2 = \theta_2^{-1}$ with larger parameter update gains}
    \label{fig:elliptical_parameter_estimates_faster_estimation}
\end{figure}

\subsection{Hilbert trajectory}
Next, the bicopter is commanded to follow a nonsmooth trajectory constructed using a second-order Hilbert curve.
The trajectory is constructed using the algorithm described in Appendix A of \cite{spencer2022adaptive} with
a maximum velocity $v_{\rm max} = 1 \ \rm m/s$ and 
a maximum acceleration $a_{\rm max} = 1 \ \rm m/s^2.$
Figure \ref{fig:Hilbert trajectory} shows the trajectory-tracking response of the bicopter, where the desired trajectory is shown in black dashes, and the output trajectory response is shown in blue. 
Figure \ref{fig:States Hilbert} shows the position $r_1$ and $r_2$ response, the roll angle $\phi$ response of the bicopter and the control $u$ generated by the adaptive backstepping controller \eqref{eq:u}. 
%

\begin{figure}[!ht]
    \centering
    \includegraphics[width=\linewidth]{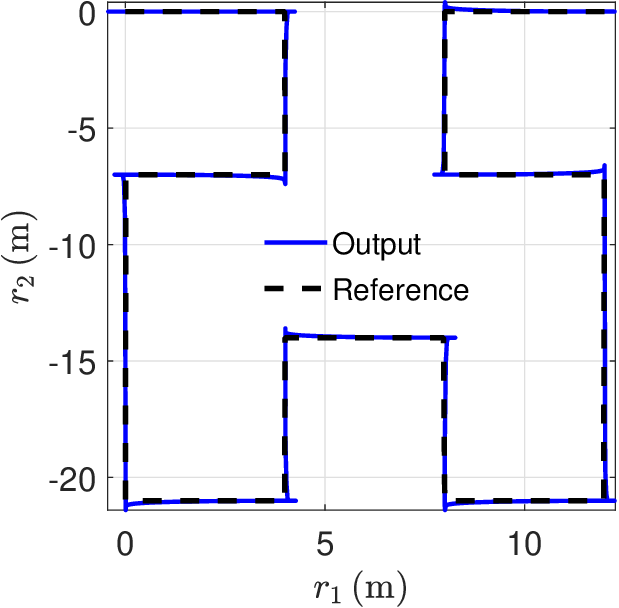}
    \caption{\textbf{Hilbert trajectory}. Tracking response of the bicopter with adaptive backstepping controller. Note that the output trajectory is in solid blue, and the desired trajectory is in dashed black.}
    \label{fig:Hilbert trajectory}
\end{figure}

\begin{figure}[!ht]
    \centering
    \includegraphics[width=\linewidth]{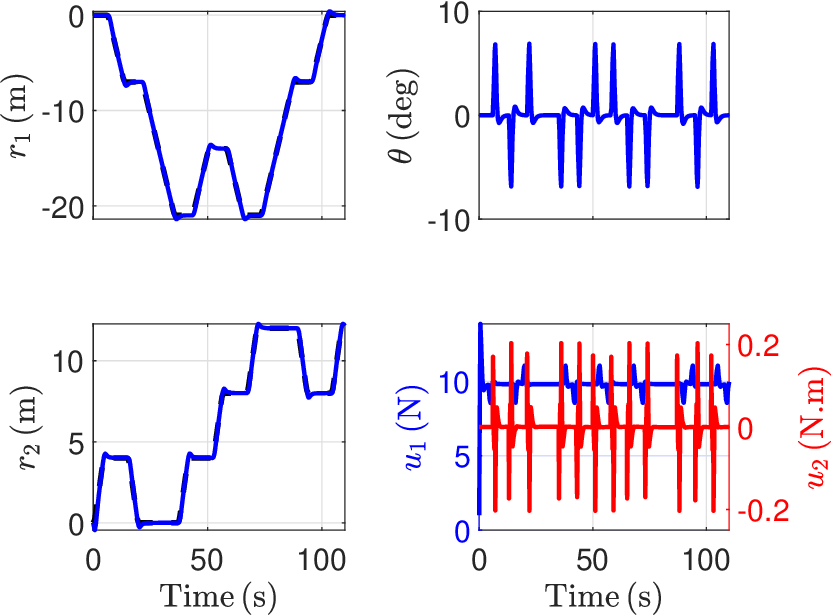}
    \caption{\textbf{Hilbert trajectory}. Positions $(r_1, r_2)$ and roll angle $\phi$ response of the bicopter with adaptive backstepping controller. Note that the references $r_{\rmd1},r_{\rmd2}$ is shown in black dashed lines}
    \label{fig:States Hilbert}
\end{figure}


Figure \ref{fig:Parameter_estimations_Hilbert} shows the estimates $\hat \theta_1,\hat \vartheta_1,$ and $\hat{p}_1$ of $\theta_1^{-1}$ and the estimate $\hat{p}_2$ of $\theta_2^{-1}.$
As in the previous case, note that the estimates do not converge to their actual values. 
%

\begin{figure}[!ht]
    \centering
    \includegraphics[width=\linewidth]{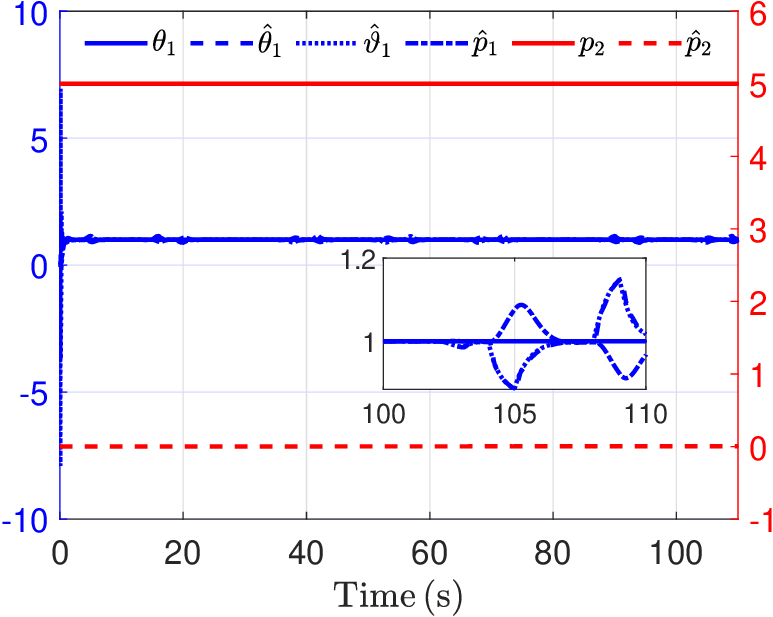}
    \caption{\textbf{Hilbert trajectory}. 
    Estimates of $\theta_1$, $\theta_1^{-1}$, and $\theta_2^{-1}$ obtained with adaption laws 
    \eqref{eq:theta_1_hat_dot}, 
    \eqref{eq:vartheta_1_hat_dot},
    \eqref{eq:p1_hat_dot}, and
    \eqref{eq:p2_hat_dot}.}
    \label{fig:Parameter_estimations_Hilbert}
\end{figure}


Next, to improve the tracking performance, we increase the estimator gains by a factor of 10; that is, the gains are $\gamma_1 = \gamma_2 =  \alpha_1 = \alpha_2 = 1.$  
Figures \ref{fig:Hilbert_trajectory_faster_estimation}, 
\ref{fig:Hilbert_trajectory_states_faster_estimation}, 
and 
\ref{fig:Hilbert_parameter_estimates_faster_estimation} show the trajectory tracking response, input signals, and the parameter estimates with the larger estimator gains. 
Note that, although the trajectory tracking error improves, the required control effort also increases.

\begin{figure}[!ht]
    \centering
    \includegraphics[width=\linewidth]{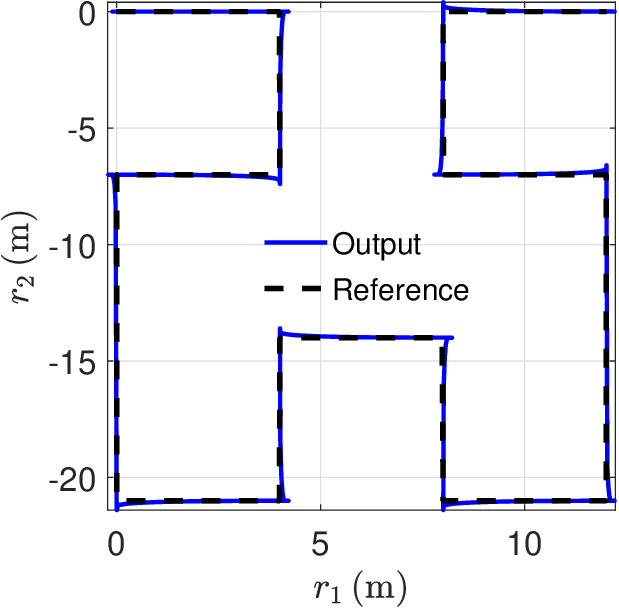}
    \caption{\textbf{Hilbert trajectory}. Tracking response of the bicopter with faster parameter estimation. 
    Note that the output trajectory is shown in solid blue, and the reference trajectory is in dashed black.}
    \label{fig:Hilbert_trajectory_faster_estimation}
\end{figure}

\begin{figure}[!ht]
    \centering
    \includegraphics[width=\linewidth]{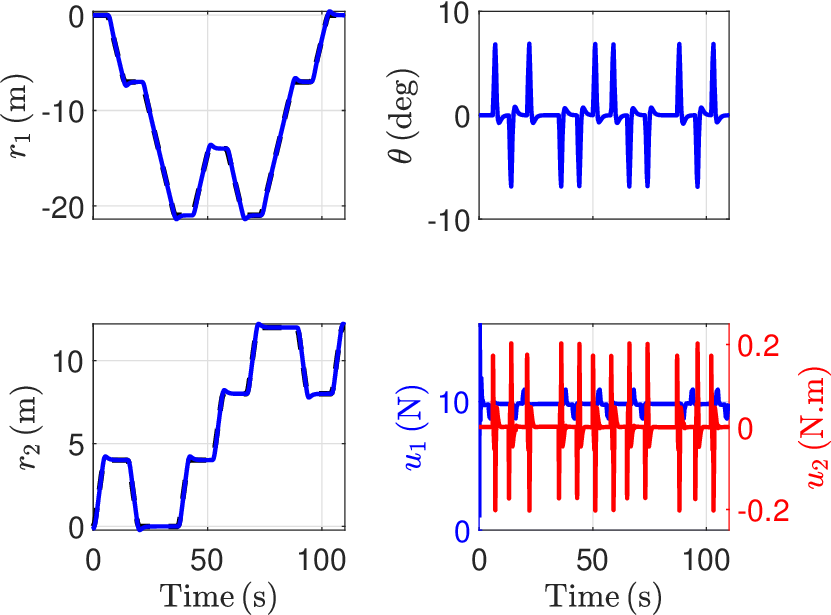}
    \caption{\textbf{Hilbert trajectory}. Position $(r_1, r_2)$ and roll angle $\phi$ response of the bicopter obtained with faster parameter estimation. Note that the references $r_{\rmd1},r_{\rmd2}$ is shown in black dashed lines}
    \label{fig:Hilbert_trajectory_states_faster_estimation}
\end{figure}

\begin{figure}[!ht]
    \centering
    \includegraphics[width=\linewidth]{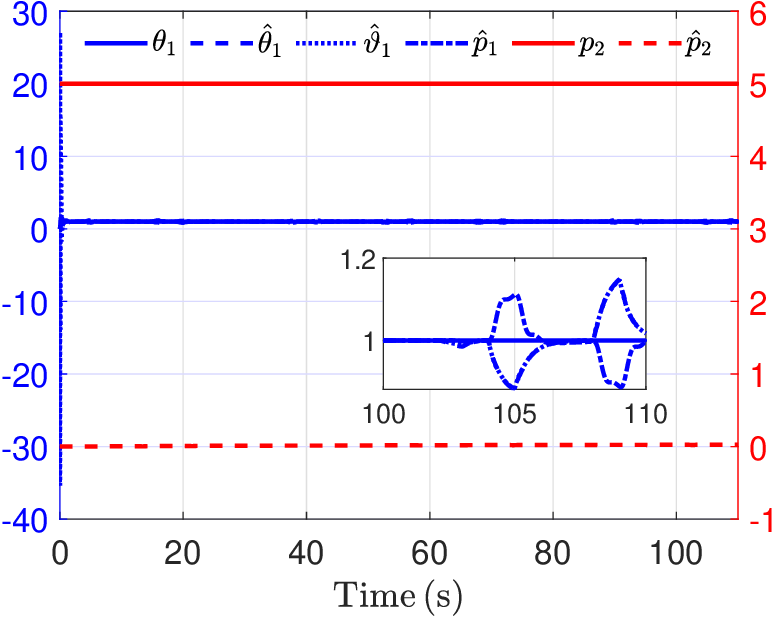}
    \caption{\textbf{Hilbert trajectory}. Estimates of $\theta_1$, $\theta_1^{-1}$, and $\theta_2$ with larger parameter update gains}
    \label{fig:Hilbert_parameter_estimates_faster_estimation}
\end{figure}

\section{Conclusions}
\label{sec:conclusions}
This paper presented an adaptive backstepping-based controller for the stabilization and tracking problem in a bicopter system. 
The adaptive controller and the parameter update laws are designed to circumvent the need to invert any system signal and only need the sign of the unknown mass and inertia, which are trivially known to be positive. 
The resulting closed-loop system is guaranteed to be stable and thus offers more flexibility in choosing gains to obtain any desired transient response. 
%
Numerical simulations show that the controller yields the desired tracking performance. 

Future work will focus on extending the proposed technique to a full multicopter model.
Furthermore, techniques to integrate state and control constraints to yield stronger performance guarantees will be investigated. 



\printbibliography
\end{document}